\documentclass [11pt]{article}
\usepackage{amsthm}
\usepackage{amssymb}
\usepackage{epsfig}

\newcounter{contador}
\newtheorem{propo}[contador]{Proposition}
\newtheorem{teo}[contador]{Theorem}
\newtheorem{lem}[contador]{Lemma}
\newtheorem{nota}[contador]{Remark}

\newtheorem{corol}[contador]{Corollary}

\newcommand{\rec}{\noindent}    
\newcommand{\dem}{\rec {\it Proof. }}  

\newcommand{\dps}{\displaystyle} 
\renewcommand{\qed}{\ \hfill\rule[-1mm]{2mm}{3.2mm}}

\newcommand{\vf}{\varphi}

\newcommand{\g}{\gamma}

\newcommand{\su}{{\mathbb S}^1} 

\newcommand{\enya}{${\rm \tilde{n}}$}

\newcommand{\lin}{{\cal L}}
\newcommand{\sg}{{\cal G}}

\newcommand{\R}{{\mathbb R}}

\newcommand{\N}{{\mathbb N}}

\newcommand{\Z}{{\mathbb Z}}
\newcommand{\M}{{\cal{M}}}

\newcommand{\U}{{\cal{U}}}

\newcommand{\CC}{{\cal{C}}}

\title{Studying discrete dynamical systems\\ through differential equations.\footnote{{\bf Acknowledgements.}  The authors are supported by MEC
through grants MTM2005-06098-C02-01 (first and second authors) and
DPI2005-08-668-C03-1 (third author). Both GSD--UAB and CoDALab
groups are supported by the Government of Catalonia through the
SGR program.}}

\author{Anna Cima$^{(1)}$, Armengol Gasull$^{(1)}$ and V\'{\i}ctor Ma\~{n}osa $^{(2)}$
  \\*[.1truecm]
{\small \textsl{$^{(1)}$ Dept. de Matem\`{a}tiques, Facultat de
Ci\`{e}ncies,}}
\\*[-.25truecm] {\small \textsl{Universitat Aut\`{o}noma de Barcelona,}}
\\*[-.25truecm] {\small \textsl{08193 Bellaterra, Barcelona, Spain}}
\\*[-.25truecm] {\small \textsl{cima@mat.uab.es, gasull@mat.uab.es}}\\
\\*[-.25truecm] {\small \textsl{$^{(2)}$ Dept. de Matem\`{a}tica Aplicada III,}}
\\*[-.25truecm] {\small \textsl{Control, Dynamics and Applications Group (CoDALab)}}
\\*[-.25truecm] {\small \textsl{Universitat Polit\`{e}cnica de Catalunya}}
\\*[-.25truecm] {\small \textsl{Colom 1, 08222 Terrassa, Spain}}
\\*[-.25truecm] {\small \textsl{victor.manosa@upc.edu}}}
\begin{document}

\maketitle
\begin{abstract}
In this paper we consider dynamical systems generated by a diffeomorphism $F$
defined on $\U$ an open subset of $\R^n,$ and give conditions over $F$ which
imply that their dynamics can be understood by studying the flow of an
associated differential equation, $\dot x=X(x),$ also defined on $\U.$ In
particular the case where $F$ has $n-1$ functionally independent first
integrals is considered. In this case  $X$ is constructed by imposing that it
shares with $F$ the same set of first integrals and that the functional
equation $\mu(F(x))=\det((DF(x))\,\mu(x),$ $x\in\U$  has some non-zero
solution, $\mu.$ Several examples for $n=2,3$ are presented, most of them
coming from several well-known difference equations.
\end{abstract}

\rec {\sl 2000 Mathematics Subject Classification:} 37C05, 37C27, 37E10, 39A20.

\rec {\sl Keywords:} Conjugation of flows, integrable vector
fields, integrable mappings, difference equations.

\section{Introduction}

In \cite{Z}, Zeeman considered the two dimensional  map
$$F(x,y)=\left(y,\frac{a+y}{x}\right)\quad {\rm with}\quad a>0,$$ which is a diffeomorphism of
$\U=\{(x,y)\in\R^2\,:\,x>0,y>0\}.$ The study of this map is motivated from the
study of the well-known Lyness difference equation $x_{k+2}=(a+x_{k+1})/x_k.$
The map $F$ has a unique fixed point in $\U$, and it has a first integral $V$,
whose level curves are topological circles surrounding the unique critical
point. By using classical tools of algebraic geometry, he proved that the
action of $F$ on each topological circle is conjugated to a rotation of the
circle, see \cite{Z} and also the paper of Bastien and Rogalsky \cite{BR}.
Zeeman conjectured that the rotation number associated to each of these
rotations, thought as a function of the energy levels of $V$ is a monotonic
function. In the end of this paper Zeeman says: ``The conjecture was initially
formulated in an attempt to give a dynamical systems proof of the result about
rotations, but it was later bypassed by finding the
 geometric proof using the cubic nature of $V.$ It might be
possible to prove the conjecture by finding a Hamiltonian system of which $F$
was the the time-1 map; this would have the advantage of placing $F$ in a
larger context." Later, in \cite{BC}, Beukers and Cushman  proved the
conjecture about monotonicity introducing the following integrable system:
\begin{equation}\label{edoBC}
\cases{\begin{array}{l}
\dot{x}=-(x+1)\left(y-\dps{\frac{x+a}{y}}\right),\\
\dot{y}=(y+1)\left(x-\dps{\frac{y+a}{x}}\right).
\end{array}
}
\end{equation}
Each energy level $\{V=h\}$  diffeomorphic to $\su$ (we  write for shortness,
$\{V=h\}\cong\su$) is a periodic orbit of the above system with period $T(h)$.
They show that on each energy level $\{V=h\}\cong\su$, there exists a number
$\tau(h)$, such that  $F$ is the  map given by the solution of system
(\ref{edoBC}) which passes through the point at a this time $\tau(h).$ They
identify the rotation number in each of these level sets as
$\rho(h)=\tau(h)/T(h)$ and in order to prove the monotonicity of $\rho(h)$ they
use the Abelian integrals--Picard--Fuchs type equation technique.

Motivated by the approach of \cite{BC}, we  follow the advice of Zeeman by
placing $F$ in a larger context. Before stating our main results we introduce
some notation and definitions.

Along this paper $\U\subset\R^n$ will be an open connected set,   $F$ will be a
diffeomorphism,  defined from $\U$ into $\U,$ and  $\dot x=X(x)$ will be a
${\cal C}^1-$differential equation also defined on  $\U.$ As usual, we denote
by $\varphi(t,p)$ its solution satisfying that $\varphi(0,p)=p,$ defined on its
maximal interval $I_p\subset\R.$  Given  $p\in \U$ we  write
$\gamma_p:=\{\varphi(t,p)\in\R^n\,,\,t\in I_p\}$ and
 ${\mathcal O}_p:=\{F^k(p)\,,\,k\in\Z\}.$

Given a diffeomorphism $F:\U\to \U,$ we will say that:

\begin{itemize}

\item the map $F$ satisfies {\it condition $X$} if there exists a ${\cal C}^1-$vector
field $X,$ also defined in $\U,$ such that for any $p\in \U,$
$$
X(F(p))=(DF(p))\,X(p).
$$

\item the map $F$ satisfies {\it condition $\mu$} if there exists a smooth map
$\mu:\U\to \R$ such that for any $p\in \U,$
$$
\mu(F(p))=\det(DF(p))\,\mu(p).
$$

\end{itemize}

Notice that condition $X$  says that the differential equation $\dot x= X(x)$
is invariant by the change of variables $u=F(x).$ On the other hand condition
$\mu$ implies several properties for the dynamical system generated by $F.$
Before listing some of them we recall a well known fact: any continuous
positive function $\nu:{\cal W}\to\R^+$ induces an absolute continuous measure
$m_\nu$ over the open set ${\cal W}\subset \U$ defined as $m_\nu({\cal B
})=\int_{\cal B} \nu(p)\,dp,$ where ${\cal B}\subset {\cal W}$ is any Lebesgue
measurable set. Recall also that it is said that $m_\nu$ is {\it an invariant
measure} for the dynamical system generated by $F$ if  for all measurable sets
$\cal B,$ $m_\nu(F^{-1}({\cal B}))=m_\nu({\cal B}).$ These are the properties:

\begin{enumerate}
\item  The set $\M^{0}:=\{p\,|\,\mu(p)=0\}$ is invariant by $F.$

\item  If  $\det(DF(p))>0$ the sets ${\cal M}^{+}:=\{p\,|\,\mu(p)>0\}$ and
${\cal M}^{-}:=\{p\,|\,\mu(p)<0\}$ are also invariant by $F$ and
$\nu(p):=\pm\frac{1}{\mu(p)}$ induces,  on each of them, an invariant
absolutely continuous measure,  $m_\nu,$ for the dynamical system generated by
$F.$

\item  If  $\det(DF(p))<0$ the sets $\M^{+}$ and $\M^{-}$ are invariant by
$F^2$ and $\nu(p)$ also induces on them, an invariant absolutely continuous
measure, $m_\nu,$ for the dynamical system generated by $F^2.$

\end{enumerate}

Property 1 is straightforward and the last two properties are
easy consequences of the Change of Variables Theorem. Notice also
that some of the above sets $\M^0,\M^+$ or $\M^-$ can be empty.
Moreover,  if we know the existence of an absolutely continuous
invariant measure for $F,$ its density gives us a way of checking
that $F$ satisfies condition $\mu.$

Observe also that if instead of condition $\mu,$ there exists a function $\mu$
such that $F$ satisfies the following condition
$$
\mu(F(p))=-\det(DF(p))\,\mu(p),
$$
then the map $F^2$ satisfies condition $\mu,$ see also Lemma \ref{eqfun}.

Finally recall that a first integral of the dynamical system
generated by a map $F:\U\rightarrow \U$ is a non-constant
$\R$-valued function $V$ which is constant on the orbits of $F.$
That is, $V(p)=V(F(p))$ for all $p\in\U.$ A set
$V_1,V_2,\ldots,V_\ell$ of first integrals of $F$ is said to be
\textsl{functionally independent} if the rank of the matrix
$(DV)(p)$ is $\ell$ for almost all $p\in \U,$ where
$V=(V_1,V_2,\ldots,V_\ell).$ Notice that for any
$\alpha\in\R^\ell$ the set $\{p\,|\, V(p)=\alpha\}$ is invariant
by $F.$

Clearly, the maximum number $\ell$ of functionally independent
first integrals for $F$ is $n.$ When $\ell=n$ it can be proved
that in most cases there exists some $k\in\N$ such that $F^k=Id,$
see \cite{CGM1.5}. The case $\ell=n-1$ is sometimes called in the
literature, {\it the integrable case,} and  will be specially
studied in this paper. The values $\alpha\in\R^{\ell}$ such that
there are points in $\{p\,|\, V(p)=\alpha\}$ where the matrix
$(DV)(p)$ has rank smaller than $\ell$ are called {\it critical
values}. The subset ${\cal C}$ of $\U$ formed by the points
belonging to level sets given by critical values will be called
{\it critical set.} It is clear that $\cal C$ is invariant by $F$
and  that
$${\cal C}=\bigcup_{\{\alpha\quad \textrm{critical}\}} \{p\,|\, V(p)=\alpha\}$$

Our first result is:

\begin{teo}\label{main1} Assume that $F:\U\to \U$ is a diffeomorphism satisfying condition $X$
and that an orbit $\gamma_p,$ solution of the differential equation $\dot
x=X(x),$ is invariant by $F.$ Then,
\begin{itemize}
\item[(a)] If $\g_p$ is a critical point of the differential equation then $p$ is a fixed point of $F.$

\item[(b)] If $\g_p$ is a periodic orbit of the differential equation then $F$
restricted to $\g_p$ is conjugated to a rotation of the circle. Moreover its
rotation number is $\tau/T,$ where $T$ is the period of $\g_p$ and $\tau$ is
defined by the equation $\varphi(\tau,p)=F(p).$

 \item[(c)] If
$\g_p$ is diffeomorphic to the real line then $F$ restricted to $\g_p$ is
conjugated to a translation of the line.
\end{itemize}
\end{teo}

Notice that in the above result, when $p$ is a critical point of $X$ it is not
necessarily a fixed point of $F.$ For instance $X(x)\equiv0$ satisfies
condition $X$ and gives no information about $F.$

For short, in the next three results, when on a regular solution $\g_p$ of
$\dot x=X(x),$ the map $F$ restricted to it satisfies one of the above last two
possibilities, namely (b) or (c), we will say that the  dynamics of $F$
restricted to $\g_p$ is {\it translation-like}.

As a consequence of Theorem \ref{main1} we prove:

\begin{teo}\label{main2}
Let $F:\U\rightarrow \U$ be a diffeomorphism having $n-1$
functionally independent first integrals $V_1,\ldots,V_{n-1}$ and
satisfying condition $\mu.$ Let $X_\mu$ be the vector field
defined as
\begin{eqnarray}
X_{\mu}(x)&=&\mu(x)\left(-\frac{\partial V_{1}(x)}{\partial x_2},\frac{\partial
V_{1}(x)}{\partial x_1}\right)\mbox{ if  } n=2,\mbox{ and}\label{casnigualados}\\
X_{\mu}(x)&=&\mu(x)\left(\nabla V_1(x)\times \nabla V_2(x)\times\cdots\times
\nabla V_{n-1}(x) \right)\mbox{ if } n>2\label{multiple}
\end{eqnarray}
where $\times$ means the cross product in $\R^n$. Then for each
regular orbit $\g_p$  of $\dot x=X_{\mu}(x),$ such that the number
of connected components of $\{x\,|\,V(x)=V(p)\}$ is $M<\infty$
there is a natural number $m,$ $1\le m \le M$ such that $\g_p$ is
invariant by $F^m$, and the  dynamics of $F^m$ restricted to
$\g_p$ is translation-like.
\end{teo}

Notice that the hypothesis that $\g_p$ is a regular orbit implies
that $p\not\in \M^0 \cup \CC.$  Observe also that the above result
forces that the periodic points of $F$ are either contained in the
set of critical points of the associated vector field $X$ or in
the periodic orbits of $X$ (which of course have to be contained
in the level sets of $V$ diffeomorphic to $\su$). Moreover they
correspond to the periodic orbits where the rotation number of $F$
is rational. In consequence these type of periodic points never
appear isolated.

We also point out that    Example 4 of Subsection \ref{exem2} shows that the
bijectivity of $F$ can not be removed from the hypotheses of the above results.

The above result can be applied to the following family of difference equations
\begin{equation}\label{eed}
x_{k+n}=\frac{R(x_{k+1},x_{k+2},\ldots, x_{k+n-1})}{x_k},
\end{equation}
where $R:{\cal U}\to \R^+,$  ${\cal U}=\{y\in \R^{n}\,|\,
y_i>0\,,\, i=1,2,\ldots, n\}$ and $n\ge 2,$ under the hypothesis
that they have $n-1$ functionally independent rational invariants.
Several examples for $n=3$ of this situation are presented in
\cite{HKY,MT}. We prove

\begin{corol}\label{coro1} Consider ${\cal U}=\{x\in \R^{n}\,|\, x_i>0\,,\,
i=1,2,\ldots, n\}.$ Let $F:\U\to \U$ be the map
$$F(x_1,x_2,\ldots,x_n)=\left(x_2,x_3,\ldots,x_{n},\frac{R(x_2,x_3,\ldots,x_n)}{x_1}\right),$$
associated to the difference equation (\ref{eed}) and assume that it has $n-1$
functionally independent rational first integrals in $\U.$ Let $\gamma_p$ be a
regular  solution of the differential equation given in
 (\ref{casnigualados}) or (\ref{multiple}),according whether $n=2$ or $n>2,$ with $\mu(x)=x_1x_2\cdots x_n.$
 Then there exists an $m\in \N,$  which is at most the number of connected components of
 $\{x\,|\, V(x)=V(p)\}$, such that the dynamics of  $F^m$ restricted to
$\g_p$ is {\it translation-like}. Moreover, when $n$ is odd $m$ has to be even.
\end{corol}

In particular, notice that the above corollary says that when $n$ is odd, all
the periodic points of $F$ of odd period have to be contained in the critical
set $\CC,$ which is an algebraic set of measure zero.

  Theorem \ref{main2} can also be
applied to the area preserving maps $F,$ giving:

\begin{corol}\label{coro2}
Let $F:\U\to \U,$ $U\subset \R^n$ be an area preserving map, {\it i.e.}
$\det(DF(x))\equiv 1$, and assume that it has $n-1$ functionally independent
rational first integrals, $V_1, V_2,\ldots, V_{n-1},$ in $\U.$ Take any smooth
function $\Phi:\R^{n-1}\to \R$ and let $\gamma_p$ be a regular solution of the
differential equation given in (\ref{casnigualados}) or (\ref{multiple}),
according whether $n=2$ or $n>2,$ with
$\mu(x)=\Phi(V_1(x),V_2(x),\ldots,V_{n-1}(x)).$ Then there exists an $m\in \N,$
which is at most the number of connected components of
 $\{x\,|\, V(x)=V(p)\}$, such that the dynamical system generated by $F^m$
restricted to $\g_p$ is {\it translation-like}.
\end{corol}

From the proof of the above results we will see  that for the maps $F$
satisfying the corresponding hypotheses and when the orbit $\g_p$ associated to
$X$ is diffeomorphic to a circle then the rotation number of $F^m,$ for some
$m\in \N,$ can be studied analytically through   $X,$ see the proof of Theorem
\ref{main1}, or \cite{BC,CGM2}. Even more, these results  open the possibility
for the numerical explorations of the rotation numbers by only using the
standard numerical methods of integration of ordinary differential equations,
see again \cite{CGM2}.

This paper is organized as follows: Section \ref{secciorelating}
is devoted to give some preliminary results. Section \ref{mr}
includes the proofs of all the results stated in this
introduction.  Section \ref{exem} deals with some concrete
applications of our results, mainly to difference equations. In
particular, Subsection \ref{exem2} deals with the case $n=2$ and
includes a revised version of the study of the planar Lyness map,
the study of the so called Gumovski--Mira-type maps and two more
examples. Subsection \ref{exem3} includes two three dimensional
maps: the Todd map, studied much extensively in \cite{CGM2}, and
an example extracted from \cite{HKY} and \cite{MT}. Finally
Section 5 contains some figures.

\section{Preliminary  results}\label{secciorelating}

Next results give some interpretations of the fact that $F$ satisfies condition
$X.$

\begin{lem} \label{fluxosconjugats}
(i) Let $F(x,y)$ be a diffeomorphism from $U$ to $U$ and assume that it
satisfies condition $X,$ {\it i.e.} $X(F(p))=(DF(p))\,X(p)$ for all $p\in U.$
Then $I_p=I_{F(p)}=I$ and for all $t\in I,$ $F(\vf(t,p))=\vf(t,F(p)).$

(ii) Conversely, if for each $p\in\U$ it holds that $F(\vf(t,p))=\vf(t,F(p))$
for $|t|$ small then $X(F(p))=(DF(p))\,X(p).$
\end{lem}

\dem (i) Consider the change of variables $u=F(x).$ Since $F$ satisfies
condition $X,$ the system $\dot x=X(x)$ becomes
$$\dot u=(DF(x))\,\dot x=(DF(x))\,X(x)=X(F(x))=X(u).$$
Since a change of coordinates gives us a conjugation of the corresponding
flows, we get $F(\vf(t,p))=\vf(t,F(p))$ whenever the equality has sense, i.e.,
for all $t\in I_p\cap I_{F(p)}.$ But since for all $t\in I_p,$ $F(\vf(t,p))$ is
well defined and it is equal to $\vf(t,F(p))$ we get that $I_p\subset
I_{F(p)}.$ On the other hand since $F$ is a homeomorphism we get that calling
$q=F(p),$ $(DF(F^{-1}(q)))\,(DF^{-1}(q))=Id$ and hence
$X(F^{-1}(q))=(DF^{-1}(q))\,X(q),$ that is, $F^{-1}$ also satisfies condition
$X.$ So, $I_{F(p)}\subset I_p$ and the equality holds.

(ii) Taking derivatives with respect to $t$ in $F(\vf(t,p))=\vf(t,F(p))$ and
substituting at $t=0$  we get the desired result. \qed

\medskip

From Lemma \ref{fluxosconjugats} we see that if $F$ is a diffeomorphism
satisfying condition $X$  then $F$ maps orbits of $X$ into orbits of $X.$ Next
proposition proves that in the case that an orbit $\gamma$ of $X$ is invariant
by $F,$ then the action of $F$ on $\gamma$ can be thought as the flow of $X$ at
a certain fixed time that only depends on $\gamma.$

\begin{propo}\label{basic}
Assume that $\gamma_p$ is invariant by $F$ and let $\tau(p)$ be defined by
$\varphi(\tau(p),p)=F(p).$ Then
$$
\varphi(\tau(p),q)=F(q)\quad\mbox{ for  all }\quad q\in\gamma_p
$$
if and only if
$$
X(F(q))=(DF)_q\,X(q)\quad\mbox{ for  all }\quad  q\in\gamma_p.
$$
\end{propo}

\dem In the following we denote $\tau=\tau(p)$. Let $q\in\gamma_p,$ that is,
$q=\varphi(t,p)$ for some $t\in I_p.$ Then if $\varphi(\tau,q)=F(q)$ for all
$q\in\gamma_p$ we have that
$$\varphi(\tau,\varphi(t,p))=F(\varphi(t,p))\quad\mbox{ for all }\quad t\in I_p.$$
Taking derivatives with respect to $t$ we get
$$X(\varphi(t+\tau,p))=(DF(\varphi(t,p))\,X(\varphi(t,p))\quad\mbox{ for  all } \quad t\in I_p,$$ which implies that
$X(F(q))=(DF(q))\,X(q)$  for all $q\in\gamma_p.$

In order to prove the converse observe that from Lemma \ref{fluxosconjugats} we
have that the function which assigns $F(\varphi(t,p))$ at each $t \in I_p\,$ is
the solution which takes the value $F(p)$ at $t=0,$ that is,
$$\varphi(t,F(p))=F(\varphi(t,p))\quad\mbox{ for  all }\quad  t\in I_p.$$ Hence, since
$\varphi(\tau,p)=F(p),$ if $q=\varphi(t,p)$ we get that
$\varphi(t,\varphi(\tau,p))=F(q)$ or, equivalently $\varphi(\tau,q)=F(q)$ for
all $q\in\gamma_p,$ as we wanted to see.\qed

\begin{nota}\label{inf}
Notice that by using the above result and Lemma \ref{fluxosconjugats} we obtain
that if $F(x,y)$ is a diffeomorphism from $U$ to $U$ satisfying condition $X,$
and it holds that  $F(p)\in\g_p$ then $I_p=I_{F(p)}=(-\infty,\infty).$
\end{nota}

It can occur that  for some point $p\in U,$ the corresponding orbit solution of
$\dot x=X(x),$ $\gamma_p$ is not invariant by $F,$ but $\gamma_p$ is contained
in a set formed by finitely many orbits which is an invariant set by $F.$ To
analyze this situation we introduce the following definition: We say that a set
of orbits of $X$, $\Gamma=\{\gamma_1,\gamma_2,\ldots,\gamma_k\}$, is an
\textsl{minimal invariant set of $F$} if $\Gamma$ is  invariant by $F$ but none
of its proper subsets is invariant by $F$.

\begin{propo}\label{minimalsets} Assume $F$ is a diffeomorphism in $\U$ and that it
satisfies  condition $X.$ Let $\Gamma$ be  a set formed by $k$ different orbits
of $X$, which is an  minimal invariant set of $F.$ Then  each $\gamma\in\Gamma$
is invariant by $F^k.$
\end{propo}

\dem From Lemma \ref{fluxosconjugats} we know that for all orbit $\gamma,$
$F(\gamma)$ is also an orbit of $X.$ Fix $\gamma_1\in\Gamma$. Observe that
$F^i(\gamma_1)\ne F^j(\gamma_1)$ for all $1\leq i< j\leq k,$ because if for
some $i$ and $j$ the equality was true,
  we would get that $\gamma_1= F^{j-i}(\gamma_1)$ and this would produce and invariant set
  formed by $j-i$ orbits,  in contradiction
with the fact that $\Gamma$ is minimal. This observation allows us to rename
the orbits in $\Gamma$  by setting $\gamma_{i+1}=F^i(\gamma),$ for
$i=1,2,\ldots,k-1$. Furthermore $F(\gamma_k)=\gamma_1$ (otherwise
$F(\gamma_k)=\gamma_i$ with $i>1$, but this would give again a contradiction).
Clearly $F(\gamma_k)=\gamma_1$ is equivalent to $F^k(\gamma_1)=\gamma_1$, hence
$F^k(\gamma_i)=F^k(F^{i-1}\gamma_1))=F^{i-1}(F^k(\gamma_1))=F^{i-1}(\gamma_1)=\gamma_i$
for all $i=1,2,\ldots,k,$ as we wanted to prove.\qed

As a corollary we get the following result:
\begin{propo}\label{periodicitat}
Assume that $F$ satisfies condition $X$ and that $p$ is a $k-$periodic point of
$F.$ Then if $X(p)\ne 0$  all the points in
$\{\gamma_p,\gamma_{F(p)},\ldots,\gamma_{F^{k-1}(p)}\}$ are $k-$periodic points
of $F.$
\end{propo}

\dem  Since $F^k(p)=p,$
$\Gamma:=\{\gamma_p,\gamma_{F(p)},\ldots,\gamma_{F^{k-1}(p)}\}$ is a minimal
invariant set. From the above proposition and Proposition \ref{basic}, all the
points in $\Gamma$ are $k$-periodic points.\qed

\bigskip

Recall that the cross product of $n-1$ vectors $W_1,\ldots,W_{n-1}$ in $\R^n,$
$n\geq 3,$ is defined  as the unique vector $W=W_1\times\cdots\times
W_{n-1}\in\R^n$ such that for all $U\in\R^n$:
$$
U\cdot W=\det\left(\begin{array}{c}
  U \\
  W_1 \\
  \vdots \\
  W_{n-1}
\end{array}
\right), $$  where the dot indicates the usual scalar product, see \cite{S}. It
is an easy exercise to check that the expression of $W$ is given by
\begin{equation}\label{defalt}
  W=W_1\times\cdots\times W_{n-1}=\sum\limits_{i=1}^n
  (-1)^{1+i}\det(M_i(W_1,\ldots,W_{n-1}))\dps{\frac{\partial}{\partial x_i}},
\end{equation}
where  $M_i(W_1,\ldots,W_{n-1})$ is the $(n-1)\times (n-1)$ matrix obtained
extracting the $i$--th column of the matrix
$$\left(\begin{array}{c}
  W_1 \\
  \vdots \\
  W_{n-1}
\end{array}
\right),$$ that is, setting
$W_k=(W_{k,1},W_{k,2},\ldots,W_{k,n}),$ the matrix given by:
$$
M_i(W_1,\ldots,W_{n-1})=\left(\begin{array}{cccccc}
  W_{1,1} & W_{1,2} &\cdots & \widehat{W_{1,i}} & \cdots & W_{1,n-1} \\
  W_{2,1} & W_{2,2} &\cdots & \widehat{W_{2,i}} & \cdots & W_{2,n-1} \\
 \vdots &  &  & \cdots &  & \vdots \\
 W_{n-1,1} & W_{n-1,2} &\cdots & \widehat{W_{n-1,i}} & \cdots & W_{n-1,n-1} \\
\end{array}
\right).
$$
We will need the following lemas:

\begin{lem}\label{lemestest}
Let $A$ be an $n\times n$ matrix, and $W_1,\ldots,W_{n-1}$ be vectors of
$\R^n$. Then
$$
  A\left(A^tW_1\times \cdots \times A^t
  W_{n-1}\right)=\det(A)\left(W_1\times\cdots\times W_{n-1}\right).
$$
\end{lem}
\dem Set $\tilde{W}=\times_{i=1}^{n} A^tW_i$, and $W=\times_{i=1}^{n} W_i$,
then, for each $U\in\R^n$ some simple computations show that
$$U\cdot
(A\tilde{W})=(A^tU)\cdot\tilde{W}=\det\left(\begin{array}{c}
  A^tU \\
  A^tW_1 \\
  \vdots \\
  A^tW_{n-1}
\end{array}
\right)=\det\left(\left(\begin{array}{c}
  U \\
  W_1 \\
  \vdots \\
  W_{n-1}
\end{array}
\right)\cdot A\right)=$$ $$=\det(A)\det\left(\begin{array}{c}
  U \\
  W_1 \\
  \vdots \\
  W_{n-1}
\end{array}
\right)=\det(A)\,(U\cdot W).$$ \qed

Notice  that in the next lemma, since the first integrals are functionally
independent, they can be labelled so that condition (\ref{nonul}) is satisfied
in a neighbourhood of most points of $\U.$

\begin{lem}\label{final}
Let $F:U\rightarrow U$ be a diffeomorphism where $U$ is an open connected set
$U\subset \R^n$ and assume that $F$ has $n-1$ functionally independent first
integrals $V_1,\ldots,V_{n-1}$. Assume that
\begin{equation}\label{nonul}
\det\left(\begin{array}{ccc}
  V_{1,1} & \cdots & V_{1,n-1} \\
  \vdots &  & \vdots \\
  V_{n-1,1} &\cdots  & V_{n-1,n-1}
\end{array}\right)(p)\ne0
\end{equation}
for all $p\in {\cal W}\subset\U,$ where $V_{i,j}=\dps{\frac{\partial
V_i}{\partial x_j}}$. Then any vector field in $\cal W$ sharing with $F$ the
same set of functionally independent first integrals writes as:
\begin{eqnarray*}
X_{\mu}(p)&=&\mu(p)\left(-V_{1,y},V_{1,x}\right)\mbox{ if  } n=2,\mbox{and}\\
X_{\mu}(p)&=&\mu(p)\left(\nabla V_1(p)\times \nabla V_2(p)\times\cdots\times
\nabla V_{n-1}(p) \right)\mbox{ if } n>2,
\end{eqnarray*}
where  $\mu:{\cal W}\to\R,$ is an arbitrary smooth function.
\end{lem}

\noindent {\it Proof.} Set $X=\sum_{i=1}^{n}
X_{i}\dps{\frac{\partial }{\partial x_i}}.$ Recall that if $X$ has
$V_k$ as a first integral then  $X(V_k)=0.$ Thus  $X$ must satisfy
the system of equations
$$  (X_1 V_{k,1}+X_2 V_{k,2}+\cdots +X_n V_{k,n})(p)=0,\mbox{ for each }
k=1,\ldots,n-1.$$ Since condition (\ref{nonul})  is satisfied in $\cal W,$
solving this system of equations by Cramer's method we have
$$
X_i(p)=\dps{\frac{X_n(p)}{\det(V(p))}}\,\det\left(\begin{array}{ccccc}
  V_{1,1} & \cdots & \overbrace{-V_{1,n}}^{i-{\rm th\, column}} & \cdots & V_{1,n-1} \\
  \vdots &  & \vdots &  & \vdots \\
  V_{n-1,1} & \cdots & -V_{n-1,n} & \cdots & V_{n-1,n-1}
\end{array}\right)(p)\, \mbox{ for }
i=1,\ldots,n-1,
$$
where $$V(p)=\left(\begin{array}{ccc}
  V_{1,1} & \cdots & V_{1,n-1} \\
  \vdots &  & \vdots \\
  V_{n-1,1} &\cdots  & V_{n-1,n-1}
\end{array}\right)(p).$$

Setting $X_n(p)=\mu(p)\det(V(p))$, we have that for each $i=1,\ldots,n-1$:
$$X_i(p)=\mu(p)\,\det\left(\begin{array}{ccccc}
  V_{1,1} & \cdots & \overbrace{-V_{1,n}}^{i-{\rm th\, column}} & \cdots & V_{1,n-1} \\
  \vdots &  & \vdots &  & \vdots \\
  V_{n-1,1} & \cdots & -V_{n-1,n} & \cdots & V_{n-1,n-1}
\end{array}\right)(p)=$$
$$=\mu(p)
(-1)^{1+i}\,\det\left(\begin{array}{cccccc}
  V_{1,1} & V_{1,2} &\cdots & \widehat{V_{1,i}} & \cdots & V_{1,n} \\
 \vdots &  &  & \cdots &  & \vdots \\
 V_{n-1,1} & V_{n-1,2} &\cdots & \widehat{V_{n-1,i}} & \cdots & V_{n-1,n} \\
\end{array}
\right)(p)=$$ $$=(-1)^{1+i}\det\left(M_i(\nabla V_1,\ldots,\nabla
V_{n-1})\right)(p).
$$ Therefore from expression (\ref{defalt}), we have that
$X(p)=\mu(p)\left(\nabla V_1\times \ldots\times \nabla V_{n-1}\right)(p),$ and
the statement is proved.\qed

By using Lemma \ref{final}, next result gives a large family  of vector fields
preserving the foliation induced by the first integrals of $F.$ It also gives
an easy way of checking whether some of  these vector fields $X$ are such that
$F$ satisfies condition $X.$

\begin{teo}\label{teoclau}
Let $F:U\rightarrow U$ be a diffeomorphism where $U$ is an open connected set
$U\subset \R^n$ and assume that $F$ has $n-1$ functionally independent first
integrals $V_1,\ldots,V_{n-1}$. Then the following statements hold:
\begin{itemize}
\item[(i)] For any smooth function $\mu:\U\to\R,$ the  family of vector fields
\begin{eqnarray*}
X_{\mu}(p)&=&\mu(p)\left(-V_{1,y},V_{1,x}\right)\mbox{ if  } n=2,\mbox{and}\\
X_{\mu}(p)&=&\mu(p)\left(\nabla V_1(p)\times \nabla V_2(p)\times\cdots\times
\nabla V_{n-1}(p) \right)\mbox{ if } n>2,
\end{eqnarray*}
shares with $F$ the same set of functionally independent first integrals.

\item[(ii)] The map $F$ satisfies condition $X_{\mu}$  if and only
if the map $F$ satisfies condition $\mu.$
\end{itemize}
\end{teo}

\noindent {\it Proof.} (i) Follows from Lemma \ref{final}.

(ii) Set $Y(p)=(\nabla V_1\times \ldots\times \nabla V_{n-1})(p)$, so that
$X(p)=\mu(p)Y(p)$. We have to prove that for all $p\in U,$
$$
{\rm (a)}\quad \mu(F(p))=\det((DF(p)))\,\mu(p) \quad\Leftrightarrow \quad{\rm
(b)}\quad X(F(p))=(DF(p)) X(p).
$$
Prior to proving this, observe that $V_i(F(p))=V(p)$ for all $i=1,\ldots,n-1$.
Then we have that $\nabla V_i(F(p))^t (DF(p))=\nabla V_i (p)^t$ and so
$$
\nabla V_i(p)= \left( \nabla V_i(F(p))^t (DF(p))   \right) ^t =(DF(p))^t \nabla
V(F(p)).
$$
Hence, we can write $Y(p)=\left((DF(p))^t\nabla V_1(F(p))\times\cdots\times
(DF)_p^t\nabla V_{n-1}(F(p))\right)$.

Using Lemma \ref{lemestest} we obtain
\begin{eqnarray*} (DF(p))
Y(p)&=&(DF)_p\left((DF)_p^t\nabla V_1(F(p))\times\cdots\times
(DF(p))^t\nabla V_{n-1}(F(p))\right)=\\
&&\det((DF(p)))\left(\nabla V_1(F(p))\times\cdots\times \nabla
V_{n-1}(F(p))\right)=\\&&\det((DF(p))) Y(F(p)).
\end{eqnarray*}

So we have proved the following identity
\begin{equation}\label{hemdeprovar} (DF(p)) Y(p)=\det((DF(p)))
Y(F(p)).
\end{equation}
Assume that condition (a) holds. Using equation (\ref{hemdeprovar}) we have
\begin{eqnarray*}
X(F(p))&=&\mu(F(p))Y(F(p))=\det((DF(p)))\mu(p)Y(F(p))=\\&&(DF(p))\mu(p)Y(p)=(DF(p))
X(p), \end{eqnarray*} hence condition  (a) implies condition (b).

Assume now that (b) holds. Then we have that $\mu(F(p))Y(F(p))=(DF(p)) \mu(p)
Y(p)$. From equation (\ref{hemdeprovar}), and since $F$ is a diffeomorphism, we
have $$ Y(p)=\det((DF(p))) (DF)^{-1}_p Y(F(p)), $$ and so
\begin{eqnarray*}
\mu(F(p))Y(F(p))&=&(DF(p))\mu(p)\det((DF(p)))(DF(p))^{-1}Y(F(p))=\\&&\mu(q)\det((DF(p)))Y(F(p)).
\end{eqnarray*}
Hence $ \mu(F(p))=\det((DF(p)))\mu(p), $ and thus condition  (b) implies
condition (a), and the proof of statement (ii) is completed.\qed

In the previous Theorem we have seen the role of condition $\mu$
for maps $F$ having $n-1$ functionally independent first
integrals. Next lemma gives some properties  that can help to find
solutions of the functional equations $
\mu(F(p))=\pm\det((DF(p)))\mu(p).$ Its proof is straightforward.

\begin{lem}\label{eqfun} Fix a diffeomorphism $F$ from $\U$ into itself, and consider
the two sets of continuous functions $\Sigma^{\pm}_F:=\{\mu:\U\to \R \,|\,
\mu(F(x))=\pm\det(DF(x))\mu(x) \}.$ Then
\begin{enumerate}
\item The spaces $\Sigma^{\pm}_F$ are  vectorial spaces.

\item If $\mu\in\Sigma^+_F$ then $\mu\in\Sigma^+_{F^k}$ for any $k\in \N.$

\item If $\mu\in\Sigma^-_F$ then $\mu\in\Sigma^+_{F^{2k}}$ for any $k\in \N.$

\item If $\mu,\nu\in\Sigma^+_F$ then $\mu^\ell\nu^{1-\ell}\in\Sigma^+_F$ for
any $\ell\in\Z.$

\item If $\mu\in\Sigma^+_F$ and $\nu\in\Sigma^-_F$ then
$\mu^\ell\nu^{1-\ell}\in\Sigma^+_{F^{2k}}$ for any $k\in\N$ and any
$\ell\in\Z.$

\item If $\mu\in\Sigma^\pm_F$ and $V$ is a first integral of $F$ then $\mu\cdot
V\in\Sigma^\pm_F.$

\item If $\det (F(x))\equiv 1$ then $\Sigma^+_F$ is the set of first integrals
of $F$ plus the constant functions.
\end{enumerate}

\end{lem}

\section{Proof of the main results}\label{mr}

\rec{\bf Proof of Theorem \ref{main1}.} $(a)$ It is trivial.

$(b)$ If $\g_p\cong \su,$ then $\g_p$ is a periodic orbit of $\dot x=X(x).$ Let
$T(p)$ be the minimal period of $\g_p.$ By proposition \ref{basic} we know that
there exists $\tau(p)$  such that $\varphi(\tau(p),q)=F(q)$ for all $q\in\g_p.$
We are going to prove that the restriction of $F$ to $\g_p$ is conjugated to a
rotation of the circle with rotation number $\rho(p)=\tau(p)/T(p).$ Let $
H:\su\longrightarrow\g_p$ and $R_\tau:\su\longrightarrow\su $ defined by $
H(e^{it})=\vf\left(\frac{T(p)}{2\pi}t,p\right)$ and $
R_\tau(e^{it})=e^{i\left(t+2\pi\frac{\tau(p)}{T(p)}\right)}.$

Then $H$ is clearly exhaustive and if $H(e^{it})=H(e^{is}),$ then
$p=\vf\left(\frac{T(p)}{2\pi}(t-s),p\right)$ and hence,
$\frac{T(p)}{2\pi}(t-s)=k\,T(p)$ for a certain $k\in\Z.$ But then $t-s=2k\pi$
and $e^{it}=e^{is}$ in $\su.$ So, it is a diffeomorphism.

To see that $H$ is the desired conjugation we have to show that
$$F\circ H=H\circ R_{\tau}.$$
As we can see, by using Lemma \ref{fluxosconjugats}, the computation is
straightforward. On one hand,
$$(F\circ H)(e^{it})\,=\,F\left(\varphi\left(\frac{T(p)}{2\pi}t,p\right)\right)=\varphi
\left(\tau(p)+\frac{T(p)}{2\pi}t,p\right),$$ and on the other one
$$(H\circ
R_{\tau})(e^{it})\,=\,H(e^{i\left(t+2\pi\frac{\tau(p)}{T(p)}\right)})\,=\,\varphi\left(\frac{T(p)}
{2\pi}\left(t+\frac{2\pi\tau(p)}{T(p)}\right),p\right)\,=\,
\varphi\left(\tau(p)+\frac{T(p)}{2\pi}t,p\right).$$

In order to prove $(c)$, assume that $\tau(p)>0$ (the case $\tau(p)<0$ follows
in a similar way). By Remark \ref{inf} we know that  $I_p=(-\infty,+\infty).$

Let
$$
H:(-\infty,+\infty)\rightarrow\g_p,\mbox{ and }
T_\tau:(-\infty,+\infty)\rightarrow (-\infty,+\infty)$$ be defined by
$H(t)=\vf(t,p),$ and $T_\tau(t)=t+\tau(p).$

Then $H(t)=H(s)$ implies $\vf(t-s,p)=p$ and since $\gamma_p$ is not a periodic
orbit we get that $t-s=0.$ On the other hand,
$$(H\circ
T_\tau)(t)=H(t+\tau(p))=\vf(t+\tau(p),p)=\vf(\tau(p),\vf(t,p))=F(\vf(t,p))=F(H(t)),$$
which proves the assertion. \qed

\bigskip

\rec{\bf Proof of Theorem \ref{main2}.} By using Theorem
\ref{teoclau} we know that the vector field $X_\mu$ introduced in
the statement of the theorem is such that $F$ satisfies condition
$X_\mu.$ Notice also that the points $p\in {\cal M}^0\cup{\cal C}$
are critical points of $X_\mu.$ Let  $p\in\U$ be a regular point
of $X_\mu.$ It is clear that $p\not\in {\cal C},$  that
$\{x\,|\,V(x)=V(p)\}$ is invariant by $F$ and has at most $M$
connected components, one of them being $\gamma_p.$ By Proposition
\ref{minimalsets} we get that for some $m\le M,$ the orbit $\g_p$
is invariant by $F^m.$ By using Theorem \ref{main1} we also know
that $F^m$ restricted to $\g_p$ is translation--like, as we wanted
to prove.\qed

\bigskip

\rec{\bf Proof of Corollary \ref{coro1}.} Note that for the map
$$F(x_1,x_2,\ldots,x_n)=\left(x_2,x_3,\ldots,x_{n},\frac{R(x_2,x_3,\ldots,x_n)}{x_1}\right),$$
the following equality holds
$$\mu(F(x))=(-1)^n\det(DF(x))\mu(x),$$
where $\mu(x)=x_1x_2\cdots x_n.$ Hence when $n$ is even,  $F$ satisfies
condition $\mu.$ On the other hand, when $n$ is odd, by using Lemma
\ref{eqfun}, we know that $F^2$ satisfies condition $\mu.$ In both cases,
 by taking $p\in\U$ a regular point for $X_\mu,$ the fact that the first
 integrals of $F$ are rational forces the finiteness of connected components of $\{x\,|\, V(x)=V(p)\}.$
Hence by using Theorem \ref{main2} the result follows.
 \qed

\bigskip

\rec{\bf Proof of Corollary \ref{coro2}.} The proof of this result
is similar to the proof of the previous corollary. The difference
is that since $F$ is a jacobian map $\mu(F(p))=\mu(p)$, so we can
take as $\mu$ any first integral of $F.$ Hence we can consider
$\Phi(V_1(x),V_2(x),\ldots,V_{n-1}(x)),$  for any  smooth function
$\Phi:\R^{n-1}\to \R.$\qed

\section{Examples}\label{exem}

\subsection{ Two dimensional examples}\label{exem2}

\noindent{\bf  EXAMPLE 1 (Lyness map).} As a first example we will see that our
approach can be applied to study the Lyness recurrence. Recall that its
associated map is
$$F(x,y)=\left(y,\frac{a+y}{x}\right)\quad {\rm with}\quad a>0,$$ defined on
$\U=\{(x,y)\in\R^2\,:\,x>0,y>0\}.$ Let us check that it is in the hypotheses of
Corollary \ref{coro1}.

It is clear that $F$ is a diffeomorphism from $\U$ to $\U$ with inverse
$F^{-1}(u,v)=((a+u)/v,u).$ It is also known that $F$ has the function
$$V(x,y)=\frac{(x+1)(y+1)(x+y+a)}{xy}$$
as a first integral, that the point $p_c=(x_c,x_c)$ with
$x_c=\frac{1+\sqrt{1+4a}}{2}$ is a fixed point of $F$ and  that the level
curves $\{V(x,y)=\alpha\}\cap \U$ are diffeomorphic to circles for all
$\alpha>V(x_c,x_c).$ Since $n=2,$ consider the vector field given in Corollary
\ref{coro1}, $X=xy\left(-\frac{\partial V}{\partial y},\frac{\partial
V}{\partial x}\right).$ We obtain
$$
X(x,y)= -(x+1)\left(y-\dps{\frac{x+a}{y}}\right)\dps{\frac{\partial}{\partial x
}}+ (y+1)\left(x-\dps{\frac{y+a}{x}}\right)\dps{\frac{\partial}{\partial y }}.
$$
Of course, it coincides with the one introduced by Beukers and Cushman in
\cite{BC}, see (\ref{edoBC}), because as we have already said that paper was
our motivating example. Since $\nabla V\ne0$ in $\U\setminus\{p_c\}$ we get
that all the orbits $\g_p$ of the differential equation associated to the above
vector field, except $p_c$ of course, are periodic orbits. Thus Corollary
\ref{coro1} implies that $F$ restricted to each of these orbits is conjugated
to a rotation.

\bigskip

\noindent {\bf EXAMPLE 2 (Gumovski--Mira-type maps)}. Consider the family of
jacobian maps
$$
F_{ \{A,B,C\}}(x,y)=\left(y,-x+\frac{By+C}{y^2+ A}\right),
$$
where $A,B$ and $C$ are real constants. It is easy to check that the family of
functions
$$
V_{ \{A,B,C\}}(x,y)=x^2y^2+ A(x^2+y^2)-B xy-C(x+y),
$$ are invariant under the action of $F_{ \{A,B,C\}}$.
By Corollary \ref{coro2}, we can consider   the vector field:
$$X_{ \{A,B,C\}}(x,y)=V_{ \{A,B,C\}}(x,y)\dps{\left(-\frac{\partial
V_{ \{A,B,C\}}}{\partial y},\frac{\partial V_{ \{A,B,C\}}}{\partial
x}\right)}.$$

A particular subfamily of the above type of maps, are the
Gumovski--Mira maps (\cite{GM}):
$$
F_{\{1,\beta,\alpha\}}(x,y)=\left(y,-x+\frac{\alpha+\beta y}{1+y^2}\right).
$$
It is easy to check that for all $\alpha,\beta$ these maps are diffeomorphisms
from $\R^2$ to $\R^2.$

From \cite{BR2} and \cite{GM}  we know when $\alpha=0$ and $\beta \in (0,2)$
the level curves of $V_{\{1,\beta,0\}}$ in $\R^2$ are closed curves surrounding
the origin, which is the unique fixed point of $F_{\{1,\beta,0\}},$ see Figure
1. These closed curves are periodic orbits of $X_{ \{1,\beta,0\}}$. It is also
easy to prove that  the only critical point of $X_{ \{1,\beta,0\}}$ is the
origin and that each level set has at most one connected component. Thus by
using Corollary \ref{coro2}, we know that on each of these closed curves
$F_{\{1,\beta,0\}},$ with $\beta \in ( 0,2),$ is conjugated to a rotation. A
detailed study of the rotation number can be found in \cite{BR2}.

When $\beta>2$ the phase portrait of the  vector field $X_{ \{1,\beta,0\}}$ is
given by two centers surrounded by a polycycle composed by a saddle with a
couple of homoclinic trajectories surrounding both centers. The rest of the
orbits are given by closed curves surrounding the polycycle, see Figure 2.

Both centers and the saddle point are fixed points of $F$. By studying the
number of connected components of the level sets of $F_{\{1,\beta,0\}}$ and by
continuity arguments  we have  that all the orbits of $X_{ \{1,\beta,0\}}$ are
invariant by $F_{\{1,\beta,0\}},$ so again from Corollary \ref{coro2} we know
that $F_{\{1,\beta,0\}}$ is conjugated to a rotation in any of the above
mentioned closed curves, and conjugated to a translation in any of the
separatrices of the polycycle.

A different situation appears for this other map studied by Gumovski and Mira
in \cite{GM}:
$$
F_{\{-a^2,-2,0\}}(x,y)=\left(y,-x+\frac{-2 y}{y^2-a^2}\right),\mbox{ with }
a\neq 1.$$ In this case the map is not defined in the whole $\R^2$. We are
going to apply again Corollary \ref{coro2}, but this time  to a given open
invariant subset of $\R^2,$ where the restriction of $F$  is a diffeomorphism.
 The level curves of the invariant $V_{\{-a^2,-2,0\}}$ have been studied in
 \cite{GM}. The results are showed in Figure 3. We can see that the origin is a center such
that its period annulus is  the open bounded set delimited by two heteroclinic
orbits which joint the two saddles.

Let $U$ be the period annulus of the origin, i. e.,
$$U=\left\{(x,y)\in\R^2:-\sqrt{a^2-1}<x<\sqrt{a^2-1},\quad\frac{-ax+a^2-1}{x-a}<y<\frac{ax+a^2-1}{x-a}\right\}.$$
Observe  $F_{\{-a^2,-2,0\}}$ is bijective on $U$ (thus a diffeomorphism). This
is  is because it is globally injective and on $U$ it has a well defined
inverse,
$F_{\{-a^2,-2,0\}}^{-1}(u,v)=\left(\frac{-2u+a^2v-u^2v}{u^2-a^2},u\right).$
However $F_{\{-a^2,-2,0\}}$ is not a diffeomorphism on $\R^2\setminus U$.

The  origin is again a fixed point of $F_{\{-a^2,-2,0\}},$ and by
continuity arguments it can be seen that the connected components
of the  closed level sets of $V_{\{-a^2,-2,0\}}$ which contain the
closed curves are invariant by $F_{\{-a^2,-2,0\}}.$ So by
Corollary \ref{coro2}, on each closed curve on the whole period
annulus U, $F_{\{-a^2,-2,0\}}$ is conjugated to a rotation of the
circle.

We like to point out that since in $\R^2\setminus U$ the map is no
more a diffeomorphism, we cannot apply our results in order to
characterize the dynamics of the map. Moreover, under the light of
Example 4 below, there is no reason to believe that each one of
the connected components of the level sets of $V_{\{-a^2,-2,0\}},$
formed by two solutions of the differential equation associated
to $X_{\{-a^2,-2,0\}},$ is invariant by $F^m$ for any $m\ge1.$

\bigskip

 \noindent {\bf EXAMPLE 3.} Consider the recurrence  $$ x_{k+2}=\frac{ax_{k+1}+b}{(c x_k+d)x_k},\quad{\rm with}\quad
a,b,c,d>0,
$$
studied in \cite{K}. It has an associated map
$$
F(x,y)=\left(y,\frac{ay+b}{(c y+d)x}\right),
$$
which is a diffeomorphism of the first quadrant $\U:=\{x>0,y>0\}$.

From \cite{K} we know that it has in $\U$ the first integral
$$
V(x,y)=\left[(d+cx)(dx+a)y^2+(a^2+bd+x^2(ac+d^2))y+(dx+a)(ax+b)\right]/(xy),
$$
and that its level curves of in $\U$ are closed curves surrounding the unique
fix point of $F$ in $\U$. By usin again Corollary \ref{coro2} we obtain that
$F$ restricted to each level set of $V$ in $\U$ is conjugated to a rotation.

\bigskip

\noindent {\bf EXAMPLE 4.} Consider the map
$$\tilde F(y,z)=\left(-\frac{1+y+z}{z},\frac{1+y+(2-a)z}{z(a+y+z)}\right),
$$ with $a\in\R.$  It is
easy to check that $H(y,z)=\frac{(1+y+z)(y+a-1)}{z}$ is a first integral of
$\tilde F,$ see \cite{CGM1}. Indeed this map is constructed in that paper by
noticing that $F^4(-1,y,z)=(-1,\tilde F(y,z)),$ where $F$ is given in formula
(\ref{lyn}) of Example 5. It is easy to verify that $\tilde F$ satisfies
condition $X$ where  $X(y,z)=X_1(y,z)\partial/\partial
y+X_2(y,z)\partial/\partial z,$ being
$$
\begin{array}{lcl}
X_1(y,z) :=-z(1+z)\frac{\partial H(y,z)}{\partial z}=(y+1)(z+1)(a-1+y)/z,\\
X_2(y,z) :=\phantom{-}z(1+z)\frac{\partial H(y,z)}{\partial y}=(z+1)(a+2y+z).
\end{array}
$$
The existence of this vector field is in fact a consequence of the study of
Example 5.

On the other hand it is proved in \cite{CGM1} that most level sets
$H_h:=\{(y,z)\,|\,H(y,z)=h\}$ are hyperbolas and that there are
many values of $h$ and $a>0$ for which the orbit of $\tilde F,$
which of course lies in a fix hyperbola, is dense over it. In
fact, in these cases, there is an appropriate compactification, so
that  $\tilde{F}$ extends to a map which is conjugated to a
irrational rotation of the circle (see Proposition 15 and the
proof of Proposition 17 of \cite{CGM1} for more details).
Therefore, although the union of the two branches of the set $H_h$
is invariant by $\tilde F,$ there is no any $m\ge1$ such that one
of the branches is invariant by $\tilde F^m.$

This example shows the importance of the hypothesis that $F$ is a
diffeomorphism in Theorems \ref{main1} and \ref{main2}.

\subsection{ Three dimensional examples}\label{exem3}

\noindent {\bf EXAMPLE 5 (Todd recurrence).} Consider the  third order Lyness
recurrence
$$
x_{k+3}=\frac{a+x_{k+1}+x_{k+2}}{x_k},\quad{\rm with}\quad
a>0,\,x_1>0,\,\,x_2>0,\,\,x_3>0,
$$
also called Todd recurrence, see \cite{CGM1, CGM2, Grove, KLR}. It has the
associated map
\begin{equation}\label{lyn}
F(x,y,z)=\left(y,z,\dps{\frac{a+y+z}{x}}\right),
\end{equation}
defined in $\U=\{(x,y,z)\in\R^3\,:\,x>0\,,\,y>0\,,\,z>0\}.$ This map is a
diffeomorphism from $\U$ into itself. From \cite{CGM2,G}, it is known that it
has a couple of functionally independent first integrals of $F$, given by:
$$\begin{array}{l}
V_1(x,y,z)=\dps{\frac { \left( x+1 \right)  \left( y+1 \right)
\left(
z+1 \right)  \left( a+x+y+z \right) }{xyz}},\\
V_2(x,y,z)=\dps{\frac{(1+y+z)(1+x+y)(a+x+y+z+xz)}{xyz}}.
\end{array}$$

In order to use  Corollary \ref{coro1} we have to consider the vector field
$$X(x,y,z)=xyz\cdot(\nabla V_1(x,y,z)\times \nabla V_2(x,y,z))=\sum_{i=1}^3 X_i(x,y,z)
\frac{\partial}{\partial x_i},$$ defined by:
$$
\begin{array}{lcl}
X_1(x,y,z) :=(x+1)(1+y+z)(yz-x-y-a)G(x,y,z)/(xy^2z^2),\\
X_2(x,y,z) :=(y+1)(z-x)(a+x+y+z+xz)G(x,y,z)/(x^2yz^2),\\
X_3(x,y,z) :=(z+1)(1+x+y)(y+z+a-xy)G(x,y,z)/(x^2y^2z),
\end{array}
$$
where $G(x,y,z)=-{y}^{3}- \left( x+a+1+z \right) {y}^{2}- \left( a+x+z \right)
y+{x}^{2 }{z}^{2}+xz+{x}^{2}z+x{z}^{2} .$ It is not difficult to check that
their critical points in $\U$ are either the points of the surface ${\cal
G}:=\{(x,y,z)\,|\, G(x,y,z)=0\}$ or the points of the line: ${\cal
L}:=\{(x,(x+a)/(x-1),x)\,|\,x>1\}.$ The line $\cal L$ is formed, precisely, by
the fix point of $F$ and the 2-periodic points of $F.$ Moreover in \cite{CGM2}
it is proved that the level curves of $(V_1,V_2)$ in $\U\setminus({\cal L}\cup
{\cal G})$ are formed by two disjoints circles. Thus if for each $p\in
\U\setminus({\cal L}\cup {\cal G})$ we consider the connected component of the
level set of $(V_1,V_2)$ which passes through it, which is diffeomorphic to a
circle, we obtain that $F^2$ restricted to it is conjugated to a rotation.
Further discussion on the rotation numbers of $F^2$, and its periods can be
found in \cite{CGM2}.

As we have seen in the proof of Corollary \ref{coro1}, the map $F^2$ satisfies
condition $\mu$ with $\mu(x,y,z)=xyz.$ In \cite{CGM2} it has been proved that
$F$ satisfies as well condition $\hat\mu$ being $\hat\mu(x,y,z)=G(x,y,z)$ the
function given above. By using Lemma \ref{eqfun} we obtain that $F^2$ satisfies
condition $\tilde \mu$ with $\tilde{\mu}=(x^2y^2z^2)/G(x,y,z).$ This is
precisely the factor for $\nabla V_1\times \nabla V_2$ used in \cite{CGM2} to
study the dynamics given by $F.$ Notice that with this new  $\tilde{\mu}$ the
 vector field is $\tilde{X}=\tilde{\mu}\cdot(\nabla V_1\times
\nabla V_2)=\sum_{i=1}^3 \tilde{X}_i \partial/\partial x_i,$ where
$$
\begin{array}{lcl}
\tilde{X}_1(x,y,z) :=(x+1)(1+y+z)(a+x+y-yz)/(yz),\\
\tilde{X}_2(x,y,z) :=(y+1)(x-z)(a+x+y+z+xz)/(xz),\\
\tilde{X}_3(x,y,z) :=(z+1)(1+x+y)(xy-y-a-z)/(xy).
\end{array}
$$
Its only critical points in $\U$ are precisely the points of the
line $\cal L$ and it allows to study the dynamics given by $F$ in
the whole $\U$ also including the surface $\cal G.$ Hence $F^2,$
and indeed $F,$ is conjugated to a rotation of the circle  on the
periodic orbits of $\tilde X$ that foliate $\cal G.$

\bigskip

\noindent {EXAMPLE 6.} In \cite{HKY}, nine third--order difference equations,
 with a couple of two functionally independent
first integrals, have been introduced. In this example we consider the
diffeomorphism
$$F(x,y,z)=\left(y,z,\displaystyle{\frac{(y+1)(z+1)}{x(1+y+z)}}\right),$$
defined in $\U=\{(x,y,z)\in\R^3\,:\,x>0\,,\,y>0\,,\,z>0\}.$  It corresponds to
a very particular case of  equation (Y$_1$) of \cite{HKY}. From \cite{RQ}, we
know that it has the following $2$-- first integrals (that is first integrals
of $F^2$), given by
$$
I_1(x,y,z)=\displaystyle{\frac{1+x+y+z+xy+yz+xyz}{xz}},\mbox{ and }
I_2(x,y,z)=\displaystyle{\frac{1+x+z+xy+xz+yz}{y}}.
$$ Consider the new two first integrals of $F$,
$$
V_1=I_1+I_2 \quad \mbox{ and }\quad V_2=I_1 I_2,
$$
see also \cite{RQ}. It is not difficult to check that the vector field that
appears in Corollary \ref{coro2}, $X(x,y,z)=xyz \cdot(\nabla V_1(x,y,z)\times
\nabla V_2(x,y,z)),$ has the surface
$$
\sg=\{G(x,y,z):=-(1+x)(1+z)y^2+[x^2z+(z^2-1)x-(1+z)]y+xz(x+1)(z+1)=0\}
$$
and the line $\lin=\{x=z; xy^2+(x^2-1)-1-x=0\},$ full of critical points.
Moreover, as in the previous example, the line $\lin$ is filled by the fix
points of $F^2,$ and $F^2$ satisfies condition $\hat\mu$ with
$\hat\mu(x,y,z)=G(x,y,z).$ At this point we can argue similarly than in the
previous case to describe the dynamics of $F^2.$

\section{Figures}

\centerline{\includegraphics[scale=0.70]{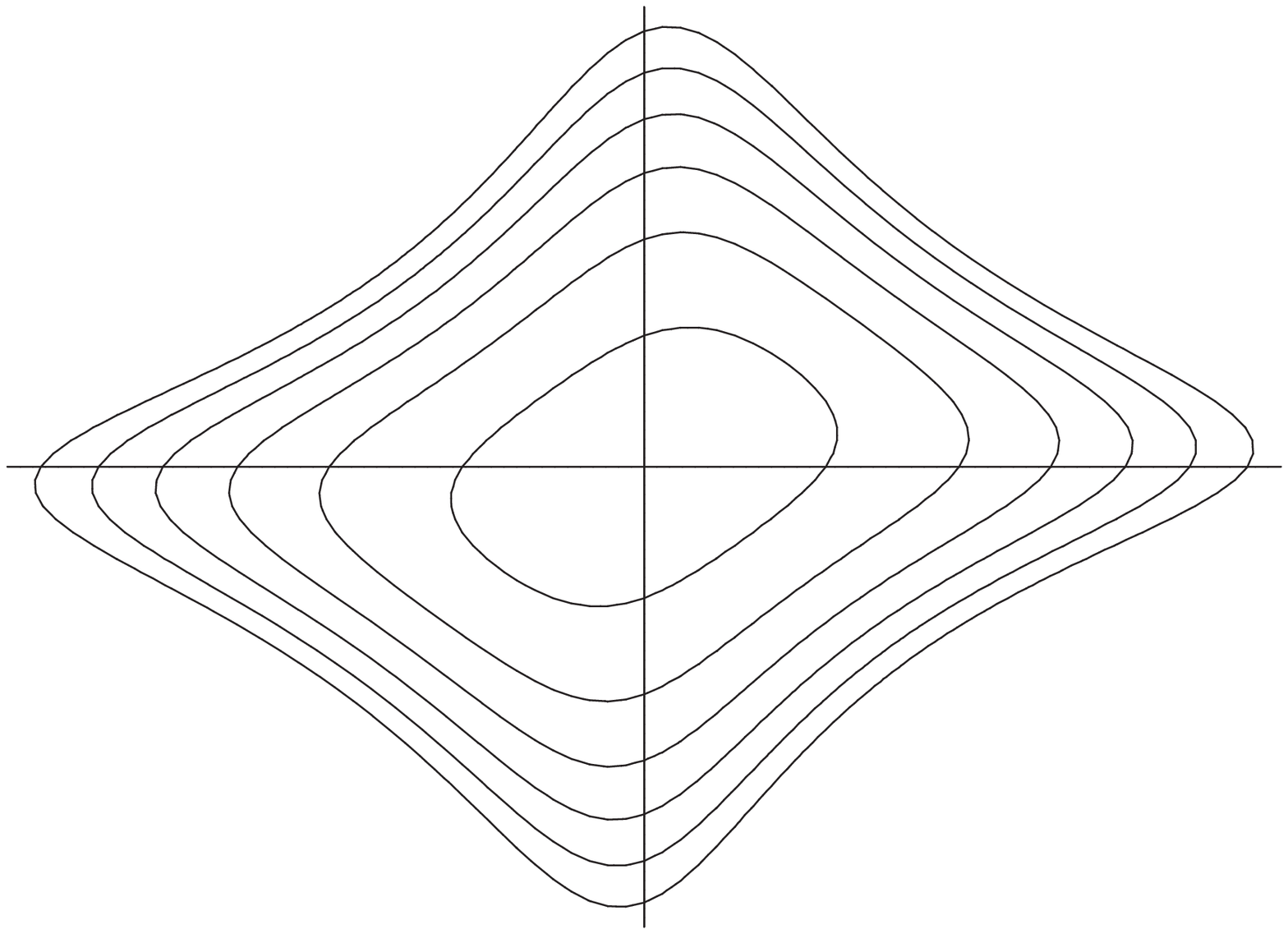}}
\begin{center}
Figure 1: Level curves of the invariant $V_{\{1,\beta,0\}}$
associated to the Gumovski--Mira map $F_{\{1,\beta,0\}}$ for
$\beta \in ( 0,2)$.
\end{center}

\centerline{\includegraphics[scale=0.70]{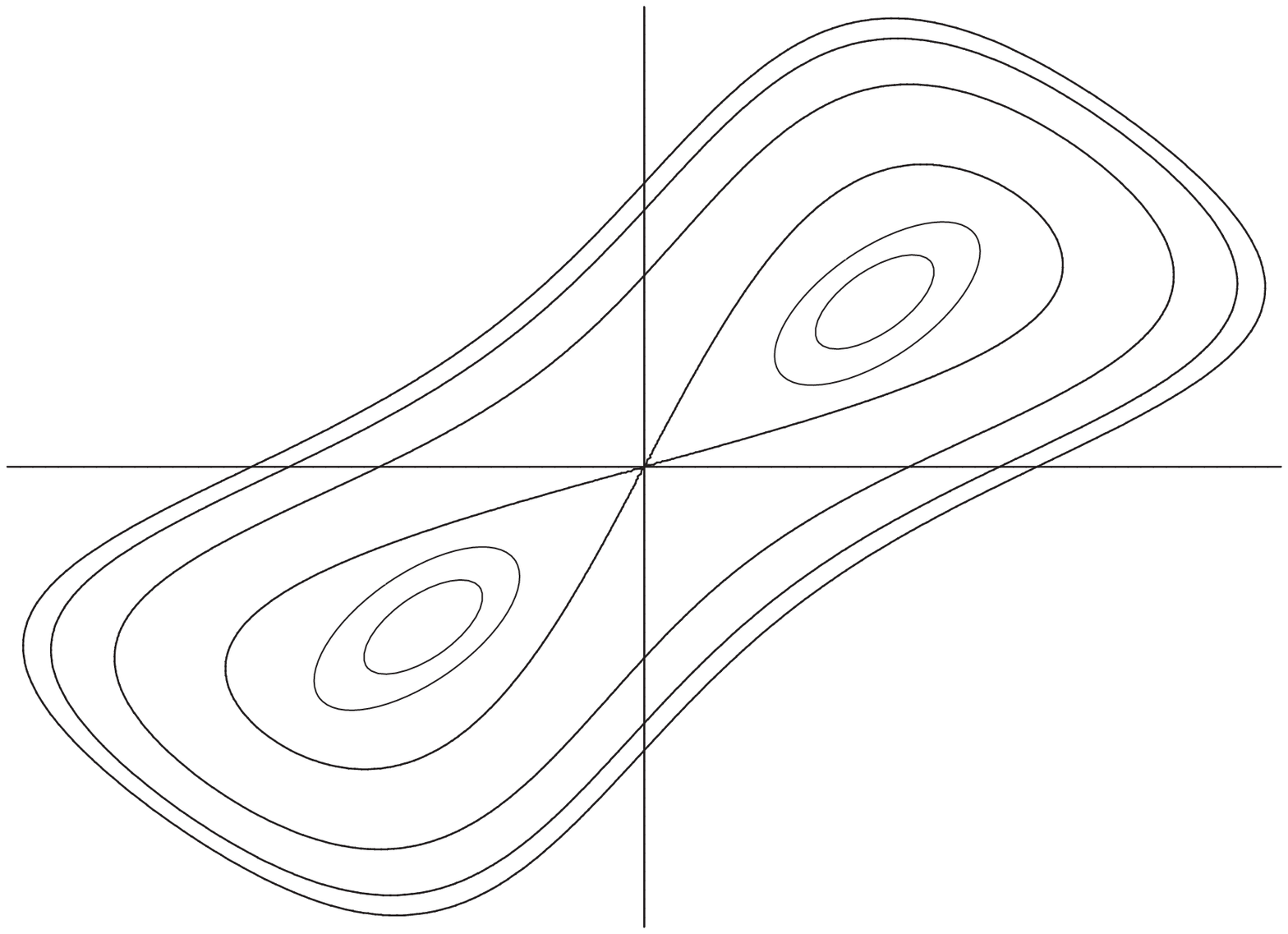}}

\begin{center}
Figure 2: Level curves of the invariant $V_{\{1,\beta,0\}}$
associated to the Gumovski--Mira map $F_{\{1,\beta,0\}}$ for
$\beta>2$.
\end{center}

\centerline{\includegraphics[scale=0.70]{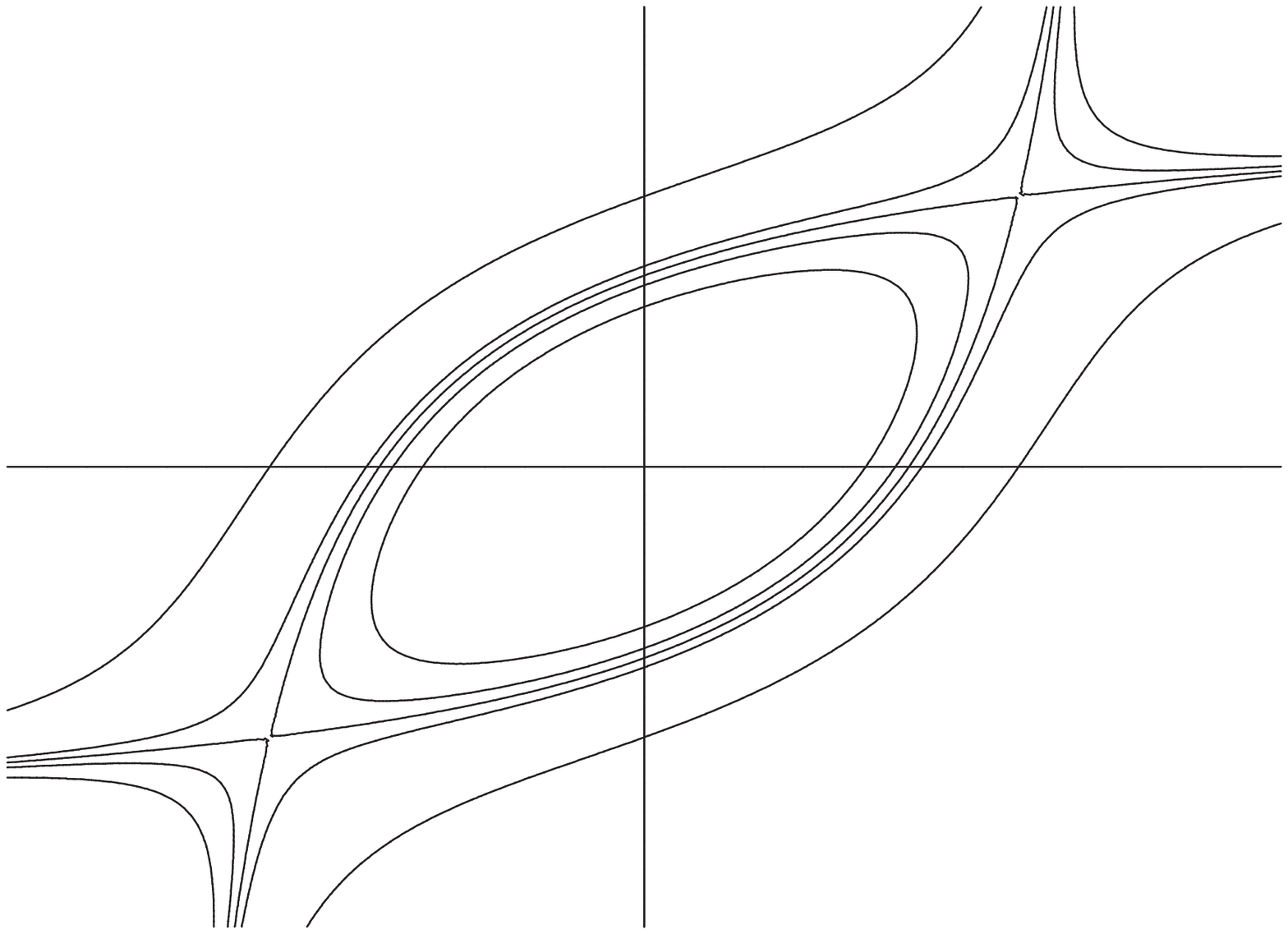}}

\begin{center}
Figure 3: Level curves of the invariant $V_{\{-a^2,-2,0\}}$
associated to the Gumovski--Mira map $F_{\{-a^2,-2,0\}}$ for
$|a|>1$. \end{center}

\end{document}